\documentclass[11pt]{article}
\usepackage{amsfonts}
\usepackage{color}
\usepackage{graphics}

\usepackage{cite}
\usepackage{latexsym}
\usepackage{amsmath}
\usepackage{amssymb}
\usepackage[dvips]{epsfig}
\usepackage{amscd}

 \usepackage{color}
 \usepackage{indentfirst}
\usepackage{ntheorem}
\theoremstyle{plain}
\theoremseparator{.}
\newtheorem{theorem}{Theorem}[section]
\newtheorem{remark}{Remark}[section]

\newtheorem{lemma}[theorem]{Lemma}

\newcommand\thmref[1]{Theorem~\ref{#1}}
\newcommand\lemref[1]{Lemma~\ref{#1}}

\newcommand{\nb}{{\bold b}}
\newcommand{\nw}{{\bold w}}
\newcommand{\ti}{\tilde}

\newcommand{\lm}{\lambda}

\def\pf{{\it Proof.}  }

\newcommand{\pa}{\partial}
\newcommand{\thatsall}{\hfill$\Box$}
\newcommand{\bi}{\bibitem}

\newcommand{\bt}{\begin{theorem}}
\newcommand{\bl}{\begin{lemma}}
\newcommand{\el}{\end{lemma}}
\newcommand{\et}{\end{theorem}}

\renewcommand{\b}{\beta  }

\newcommand{\te}{\theta}

\newcommand{\al}{\alpha}
\newcommand{\de}{\delta}
\newcommand{\ve}{\varepsilon}
\newcommand{\la}{\label}
\newcommand{\si}{\sigma}
\newcommand{\ka}{\kappa}

\newcommand{\bn}{\begin{eqnarray}}
\newcommand{\en}{\end{eqnarray}}
\newcommand{\bnn}{\begin{eqnarray*}}
\newcommand{\enn}{\end{eqnarray*}}

\newcommand{\bnnn}{\begin{eqnarray*}}
\newcommand{\ennn}{\end{eqnarray*}}
\newcommand{\ben}{\begin{enumerate}}
\newcommand{\een}{\end{enumerate}}
\newcommand{\ba}{\begin{aligned}}
\newcommand{\ea}{\end{aligned}}
\newcommand{\be}{\begin{equation}}
\newcommand{\ee}{\end{equation}}

\def\norm[#1]#2{\|#2\|_{#1}}

\def\xix{\int_0^1}

\def\xiT{\int_0^T}

\def\O{\Omega}

\makeatletter
\@addtoreset{equation}{section}
\makeatother

\makeatletter
\@addtoreset{equation}{section}
\makeatother

\title{ Large-time Behavior of     Magnetohydrodynamics  with Temperature-Dependent  Heat-Conductivity  \thanks{ Partially supported by NNSFC   11671027 and 11471321.}
}

\author{Bin Huang$^1$,  Xiaoding Shi$^1$, Ying Sun$^2$ \thanks{Email addresses: abinhuang@gmail.com, abinhuang36@163.com (B. Huang), shixd@mail.buct.edu.cn (X. Shi)   1913349041@qq.com (Y. Sun)
}  
\\[3mm]   1.  Department of Mathematics, College of Mathematics and Physics, \\Beijing University of Chemical Technology, \\ Beijing  100029, P. R. China  
\\[3mm]   2.  School of  Mathematical Sciences,   Xiamen University, \\ Xiamen 361005, P. R. China }
\date{ }

\begin{document}
\maketitle

\begin{abstract}
We study the large-time behavior of strong solutions to   the equations of a planar magnetohydrodynamic compressible flow  with  the heat conductivity    proportional to a nonnegative power of  the temperature.
Both the specific volume and the temperature are  proved to be  bounded  from below and above independently of time. Moreover, it is shown that the global strong solution is
nonlinearly exponentially stable   as time tends to infinity. 
Our result can be regarded as a natural generalization of the previous ones   for the compressible Navier-Stokes system to MHD system with either constant heat-conductivity or  nonlinear  and  temperature-dependent     heat-conductivity.   \end{abstract}

$\mathbf{Keywords.}$ Magnetohydrodynamics, Temperature-dependent     heat-conductivity,  Strong solutions, Large-time behavior  

{\bf Math Subject Classification:} 35Q35; 76N10.
%\newpage

\section{Introduction}

The governing equations of a planar magnetohydrodynamic compressible flow written
  in the Lagrange variables read as follows:
  \be\la{1.1}
v_t=u_{x},
\ee
\be\la{1.2}
u_{t}+(P+\frac12 |\nb|^2)_{x}=\left(\mu\frac{u_{x}}{v}\right)_{x},\ee
 \be\la{1.3}\nw_t-\nb_x=\left(\lm\frac{\nw_{x}}{v}\right)_{x}, \ee \be\la{1.4} (v\nb)_t-\nw_x =\left(\nu\frac{\nb_{x}}{v}\right)_{x}, \ee
\be\la{1.5}\ba&\left(e+\frac{u^2+|\nw|^2+v|\nb|^2}{2}\right)_{t}+ \left(
u\left(P+\frac12|\nb|^2\right)-\nw\cdot\nb\right)_{x}\\&\quad=\left(\kappa\frac{\theta_{x}}{v}
+\mu\frac{uu_x}{v}+\lm\frac{\nw\cdot\nw_x}{v}+\nu\frac{\nb\cdot\nb_x}{v}\right)_x,
 \ea
\ee
where $t>0$ is time, $x\in \Omega=(0,1)$ denotes the
Lagrange mass coordinate,  and the unknown functions $v>0, u, \nw\in \mathbb{R}^2, \nb \in   \mathbb{R}^2, e>0, \theta>0,$ and $P$ are,  respectively, the specific volume of the gas, longitudinal velocity, transverse velocity, transverse magnetic field, internal energy,  absolute temperature and  pressure. $\mu$ and $\lm$ are the viscosity of the flow,
$\nu$ is the magnetic diffusivity of the magnetic field, and $\ka$ is the heat conductivity.

In this paper, we
consider a perfect gas for magnetohydrodynamic flow, that is, $P$ and $e$ satisfy \be \la{1.6}   P =R \theta/{v},\quad e=c_v\theta +\mbox{const},
\ee
where  both specific gas constant  $R$ and   heat
capacity at constant volume $c_v $ are   positive constants.
We also assume that $\mu,  \lm, $ and $\nu$ are positive constants, and  $ \ka$ satisfies \be\la{1.7}   \ka=\ti\ka \te^\beta, \ee
with constants $ \ti\ka>0$ and $ \beta\ge 0.$

The system \eqref{1.1}-\eqref{1.7} is supplemented
with   initial  conditions
\be\la{1.8}(v,u,\te,\nb,\nw)(x,0)=(v_0,u_0,\te_0,\nb_0,\nw_0)(x),  \quad  x\in\Omega, \ee
and boundary   ones\be\la{1.9} \left(u,\nb ,\nw,\theta_x\right)|_{\partial\Omega}=0,\ee
where  the initial data \eqref{1.8} should be compatible with the boundary conditions \eqref{1.9}.

Magnetohydrodynamics (MHD), concerning the flow of electrically
conducting fluids in the presence of magnetic fields,
covers a wide range of physical objects from liquid metals to cosmic plasmas (\cite{q1,q2,q3,q4,q5,q6,q7}).  The central point of MHD theory is that conductive fluids can support magnetic fields. % MHD is a macroscopic theory.
The partial differential equations of MHD can in principle be derived from Boltzmann's equation assuming space and time scales to be larger than all inherent scale-lengths such as the Debye length or the gyro-radii of the charged particles (\cite{q1,q3, q5,q6,q7}). In fact,
one can deduce from the Chapman-Enskog expansion for the first level of approximation in kinetic theory that the viscosity $\mu$
and heat conductivity $\ka$ are functions of temperature alone (see Chapman-Colwing \cite{cc1}).     These dependencies, especially the dependence of viscosity on temperature, brings great difficulties and challenges to mathematical analysis and numerical calculation. Thus, to study this problem, we first consider the case that the  viscosity is a positive constant and    the heat conductivity    proportional to a nonnegative power of  the temperature,  as shown as in the equation \eqref{1.7}. % The equations \eqref{1.1}--\eqref{1.5} describe the macroscopic behavior of the three-dimensional MHD flow, which is uniform in the transverse directions, with dissipative mechanisms. The assumption \eqref{1.7} is motivated by the physical facts: the heat conductivity $\ka$ and viscosity $\mu$ of compressible heat-conducting magnetic fluid vary with temperature and density under very high temperature and density environment (see \cite{be1,ze} and the references therein).

There is huge literature on the studies of the global existence and large time behavior of solutions to the  compressible Navier-Stokes system and MHD system. Indeed,
for compressible Navier-Stokes system \eqref{1.1} \eqref{1.2} \eqref{1.5} with $\nb\equiv\nw\equiv 0$,    Kazhikhov and Shelukhin
\cite{9} first obtained   the global existence of solutions
 for constant coefficients $(\beta=0)$  with large initial data. From then on, much effort has been made to generalize this approach to other cases (for $\b>0,$ see    \cite{24,28,hs1} and the reference therein). As for MHD system, the are many results concerning the global existence of solutions  with large initial data (see \cite{ka1,az1,6,7,cw1,cw2,fjn1,25,fhl1,hss1} and the references therein).      In particular, Kazhikhov \cite{ka1} (see also \cite{az1}) first for $\b=0$ and very recently  Huang-Shi-Sun  \cite{hss1} for $\b>0$  proved that
 \begin{lemma}[\cite{ka1,hss1}] \la{thm11.1} Let $\b\ge 0.$ Suppose % that \be \la{a1.10}\al\ge 0,\quad \beta> 0,\ee
 that the initial data $ ( v_0,u_0,\te_0,\nb_0,\nw_0)$   satisfies
  \be  \la{a1.11} ( v_0,\te_0)\in   H^1 (0,1),\quad (u_0,\nb_0,\nw_0)\in H^1_0 (0,1),\ee  and \be \la{a1.12}
\inf_{x\in (0,1)}v_0(x)>0, \quad \inf_{x\in (0,1)}\theta_0(x)>0. \ee
Then, the initial-boundary-value problem \eqref{1.1}--\eqref{1.9} has a unique strong solution $(v,u,\te,$\\$\nb,\nw)$ such that for each fixed $ T>0 $,
 \be
 \begin{cases}   v ,\,\theta \in L^\infty(0,T;H^1(0,1)),\quad u,\,\nb,\nw  \in L^\infty(0,T;H^1_0(0,1)),\\ v_t\in
  L^\infty(0,T;L^2(0,1))\cap L^2(0,T;H^1(0,1)), \\ u_t,\,\theta_t,\,\nb_t,\,\nw_t,\,u_{xx},\,\te_{xx},\,\nb_{xx},\,\nw_{xx} \in
  L^2((0,1)\times(0,T)),\end{cases}\ee
  and for each $(x,t)\in[0,1]\times[0,T]$
  \be\la{pv1} C^{-1}\leq v(x,t)\leq C,\quad C^{-1}\leq\te(x,t)\leq C,\ee
  where $C>0$ is a constant  depending on the   data and $T.$
  \end{lemma}

Concerning  the large-time behavior of    solutions to the compressible Navier-Stokes system \eqref{1.1} \eqref{1.2} \eqref{1.5} with $\nb\equiv\nw\equiv 0$,    Kazhikhov \cite{kaz1} (see also \cite{akm,az11,az2,az3,na1,na2,na3,ni,qy} among others) first obtained that for the case that $\b=0,$   the strong solution   is
nonlinearly exponentially stable   as time tends to infinity.  Very recently, Huang-Shi \cite{hs1} prove that the same result still  holds for the compressible Navier-Stokes system \eqref{1.1} \eqref{1.2} \eqref{1.5} with $\nb\equiv\nw\equiv 0$ for $\b>0.$    However, it seems to us that the known lower and upper bounds
of  the the specific volume $v$ and the temperature $\te$   depend on the time $T$ , see \cite{ka1,hss1}, so it is impossible to study
the large time asymptotic behavior of solutions in the setting in \cite{ka1,hss1}. %A natural question is what is the large time behavior of  A natural question arises: What is the large-time behavior  the strong solutions to the MHD system \eqref{1.1}-\eqref{1.9} for $\b\ge 0$ %It was somewhat of a surprise to us that it remains completely open  even for $\b=0.$
In fact, the main aim of this paper is to prove that the global strong solutions whose existence is guaranteed by Lemma \ref{thm11.1} are
nonlinearly exponentially stable   as time tends to infinity for $\b\ge 0.$

We now state
  our main result   as follows.

 \begin{theorem} \la{thm11} Under the conditions of Lemma \ref{thm11.1}, there exist    positive constants $C$ and $\eta_0 $ both depending only on the   data %on $R,c_v,\mu,\nu,\lm,\ti\ka, \b,\|(v_0,u_0,\theta_0,\nb_0,\nw_0)\|_{H^1(0,1)},  \inf\limits_{x\in [0,1]}v_0(x),$ and $ \inf\limits_{x\in [0,1]}\theta_0(x) $
 such that the unique strong solution $(v,u,\te,\nb,\nw)$  of the initial-boundary-value problem  \eqref{1.1}--\eqref{1.9}  obtained by Lemma \ref{thm11.1} satisfies for any $(x,t)\in (0,1)\times(0,\infty),$
   \be C^{-1}\le v(x,t)\le C,\quad C^{-1}\le\te(x,t)\le C,\ee and for any $t>0,$
    \be\la{vp1} \ba &\left\|\left(v -v_s,u, \te-\te_s,\nb,\nw  \right)  \right\|_{H^1(0,1)}  \le C e^{-\eta_0t},\ea\ee with  $$v_s=\int_0^1v_0dx,\quad \te_s=\int_0^1\left(\te_0+\frac{u_0^2 +|\nw_0|^2+v_0|\nb_0|^2}{2c_v}\right)dx.$$
  \end{theorem}

  %A few remarks are in order.

\begin{remark} Our result can be regarded as a natural generalization of   previous ones   for compressible Navier-Stokes system    with  either constant  heat conductivity   (\cite{kaz1}) or   temperature-dependent one (\cite{hs1}) to the  temperature-dependent  heat conductivity MHD system with  the  constant  heat conductivity as a special case.\end{remark}

We now make some comments on the analysis of this paper. The key step to study the large-time behavior  of  the global strong  solutions is to get the     time-independent lower and upper bounds   of both $v$ and $\theta $ (see   \eqref{cqq1},  \eqref{vg14},  \eqref{q2.49}, and \eqref{tq2}).  Compared with \cite{kaz1,ni,hs1}, the main difficulties come  from  the interaction of  the hydrodynamic and electrodynamic effects and  the degeneracy and nonlinearity of the heat conductivity.    Hence, to overcome these difficulties, some new ideas are needed. The key observations are as follows: First, after modifying  the ideas due to   \cite{kaz1,hs1}, we obtain an   explicit  expression of the specific volume $v$ (see \eqref{2.6}) which together with a lower bound of the temperature $\te$ (see \eqref{twq1}) shows that  $v$ is bounded from   below time-independently (see \eqref{cqq1}). Then, for $\b>0,$ we find that (see \eqref{u1}) \bnn \xiT\max\limits_{x\in[0,1]}\left(\te^{1/2}(x,t)- 2\right)_+^2 dt\le C ,\enn which gives  the uniform upper bound of $v .$ For $\b=0,$ it seems much more difficult to bound $v$ from above time-independently.  We  first prove a new estimate that (see \eqref{q2.45}) $$\int_0^T\int_0^1 \frac{|\nw_x|^2}{v}dxdt\le C,$$ which play an important role in the analysis. Then we  refine the strategy of Kazhikhov (\cite{ka1}), that is, we prove that the $L^\infty(0,T; L^2(0,1))$-norm of $(\ln v)_x$ can be bounded time-independently by a  log-type inequality (see \eqref{2.45}), which together with the Gronwall inequality in turn gives the uniform upper bound of $v$ (see \eqref{q2.49}). Next, for $\b>0$ and for the upper and lower bounds of $\te,$ we modify slightly the ideas due to \cite{hs1}, that is,    we prove  that the $L^\infty(0,\infty;L^p)$-norm of $\te^{-1}$ is bounded (see \eqref{ljj1}),  which  yields that    the $L^2((0,1)\times (0,T))$-norm of $\te_x$ is bounded provided $\b>1$ (see \eqref{sq1}). Finally, for $\b\in[0,1], $ we find that  the $L^2((0,1)\times (0,T))$-norm of $\te_x$ can be  bounded by  the $L^4( 0,T;L^2(0,1))$-norm of $u_x$ which plays an important role in obtaining  the uniform bound on $L^2((0,1)\times (0,T))$-norm of both $\te_x$ and $u_{xx} $ (see Lemma \ref{lemmy11}) for $\b\in [0,1]$.   The whole procedure will be carried out in the next section.

\section{ Proof of \thmref{thm11}}

Without loss of generality, we assume that $\lambda=\nu= \mu = \ti\ka=R=c_v=1 ,  $ and that
 \bnn  \int_0^1 v_0dx=1,\quad  \int_0^1\left(\te_0+\frac{u_0^2+|\nw_0|^2+v_0|\nb_0|^2}{2}\right)dx=1.\enn

%\newpage
We first state  the time-independent  lower bound  of $v$.
\begin{lemma}\la{wq20} For $\b\ge 0,$ it holds that for any $(x,t)\in[0,1]\times [0,+\infty),$
\be\ba\la{cqq1}
v(x,t)\ge C_0  ,
\ea\ee where (and in what follows)   $C_0  $ and $C$  denote  some
generic positive constants
 depending only on $\b,\|(v_0,u_0,\theta_0,\nb_0,\nw_0)\|_{H^1(0,1)},
 \inf\limits_{x\in [0,1]}v_0(x),$ and $ \inf\limits_{x\in [0,1]}\theta_0(x).$
\end{lemma}

\pf   First, it follows from  \eqref{1.1}, \eqref{1.5},   and   \eqref{1.9}  that for $t>0$
 \be\la{wq1}  \int_0^1 v(x,t)dx\equiv 1,  \quad \int_0^1\left(\te +\frac{u ^2+|\nw |^2+v |\nb |^2}{2}\right) (x,t)dx \equiv 1.\ee

  Then,    denoting
\be\la{2.9}\ba  \si\triangleq   \frac{u_x}{v}-\frac{\te}{v}-\frac12|\nb|^2, \ea\ee % which satisfies \be \la{a2.9} \si=( \ln v )_t -\frac{\te}{v}-\frac12|\nb|^2,\ee due to \eqref{1.1},
we rewrite \eqref{1.2} as
\bnn\la{2.10}\ba u_t=\si_x.\ea\enn
Integrating  this with respect to $x$  over $(0,x)$ gives
\be\la{2.11}\ba \left(\int_0^xudy\right)_t=\si-\si(0,t) , \ea\ee
 which implies
\bnn\ba\la{2.0} v\si(0,t)=v\si-v\left(\int_0^xudy\right)_t . \ea\enn
Integrating this with respect to $x$ over $(0,1),$ we obtain after using \eqref{1.9},   \eqref{wq1},    and \eqref{2.9}  that \be\ba\la{2.13}   \si(0,t) &= \int_0^1\left(  u_x-\te-\frac{v}{2}|\nb|^2 \right) dx-\left(\int_0^1v\int_0^xudydx\right)_t\\&\quad +\int_0^1u_x\int_0^xudy dx\\&=- \left(  \int_0^1v\int_0^xudydx\right)_t -\int_0^1\left(\te+ \frac{v}{2}|\nb|^2 + u^2 \right) dx.\ea\ee

 Finally, combining \eqref{2.11}, \eqref{1.1}, and \eqref{2.13} yields
\be\la{2.12}\ba v(x,t) =   D (x,t)Y (t)  \exp\left\{ \int_0^t\left(\te+  \frac{v}{2}|\nb|^2\right)v^{-1} d\tau\right\}, \ea\ee
 with  \be\ba\la{a2.7}D (x,t)=&v_0\exp\left\{    \int_{0}^x \left(u(y,t)-u_0(y)\right)dy\right\}\\&\times \exp\left\{- \int_0^1v\int_0^xudydx+ \int_0^1v_0\int_0^xu_0dydx\right\},\ea\ee
and \be\ba\la{2.8} Y (t)=\exp\left\{  - \int_0^t\int_0^1\left(u^2+ \frac{v}{2}|\nb|^2+\te\right) dxd\tau\right\} .\ea\ee  Using \eqref{2.12}, direct computation gives
\be \la{2.6}\ba v(x,t)=   D (x,t)Y (t)  \left\{1 + \int_0^t\frac{(\te+  \frac{v}{2}|\nb|^2)(x,\tau)}{  D (x,\tau)Y (\tau)} d\tau \right\}. \ea\ee

Next,  using  \eqref{1.1}-\eqref{1.4}, we rewrite the energy equation \eqref{1.5}   as
\be\la{2.4}\theta_{t}+  \frac\theta v
u_{x}= \left(\frac{\theta^\b\theta_{x}}{v}\right)_{x}+ \frac{  u_{x}^{2}+|\nw_x|^2+|\nb_x|^2}{v}.\ee    Multiplying   \eqref{1.1},
   \eqref{1.2},       \eqref{1.3}, \eqref{1.4},   and \eqref{2.4}   by $ 1- {v}^{-1} ,
u, \nw, \nb,$ and  $  1- {\theta}^{-1} $  respectively,   adding them altogether and integrating the result over ${(0,1)\times(0,T)}$, we obtain the following energy-type inequality
 \be\ba\la{2.2} &\sup_{0\le t\le T}\int_0^1\left(  \frac{u ^2+|\nw |^2+v |\nb |^2}{2} + (v-\ln
v )+ (\theta-\ln \theta )\right)dx \\&\quad+\int_0^T V(s)ds \le e_0, \ea\ee
where
\bnn\ba V(t)\triangleq\int_0^1\left(\frac{\theta^\b\theta_{x}^2}{v\te^2}+\frac{  u_{x}^{2}+|\nw_{x}|^{2}+  |\nb_{x}|^{2}}{v\te}\right)(x,t)dx. \ea\enn and
\bnn e_0\triangleq 2\int_0^1\left(  \frac{u_0^2+|\nw_0 |^2+v_0 |\nb_0 |^2}{2} + (v_0-\ln
v_0 )+ (\theta_0-\ln \theta_0 )\right)dx. \enn

Next,  applying Jensen's inequality to the convex  function $\te-\ln \te$  leads to
\bnn\ba \int_{0}^{1}  \te dx-\ln\int_{0}^{1}   \theta dx  \le \int_{0}^{1} (\te-\ln\theta ) dx ,\ea\enn
  which together with \eqref{2.2} and \eqref{wq1} leads to \be \la{qwe1} \bar \te(t)\triangleq \int_0^{1}\te(x,t)dx\in[ \al_1,1],\ee   where $0<\al_1<\al_2$  are two roots of  \bnn x-\ln x =e_0.\enn

  Next, both \eqref{wq1} and Cauchy's inequality give
  \bnn\ba
 \left|\int_0^1 v\int_0^x udydx\right|&\le \int_0^1 v\left|\int_0^x udy\right|dx\\&\le\int_0^1 v\left(\int_0^1 u^2dy\right)^{1/2}dx\\&\le C,
 \ea\enn which combined with \eqref{a2.7} shows
\be\ba \la{cq8}C^{-1}\le  D(x,t)\le  C.\ea\ee
Moreover,  one deduces from \eqref{wq1} that\bnn\ba\la{cq9} \al_1\le \int_0^1\left(u^2+\frac{v|\nb|^2}{2}+\te\right) dx\le 2,\ea\enn
which   yields that for any $0\le \tau<t<\infty, $
\be \la{cq7}e^{-2t}\le Y(t)\le 1,\quad e^{-2(t-\tau)}\le  \frac{Y(t)}{Y(\tau)}\le e^{- \al_1(t-\tau)}.\ee

Next, denoting $f_+\triangleq \max\{f,0\},$ we have
\bnn\ba  \left(\bar\te^{\frac{\b+1}{2}}( t)- \te^{\frac{\b+1}{2}}(x, t) \right)_+ &\le \int_0^1 \left|\pa_x\left(\bar\te^{\frac{\b+1}{2}}( t)- \te^{\frac{\b+1}{2}}(x, t) \right)_+\right|dx\\&\le C\left(\int_0^1 \frac{ \te^\b \te_x^2}{\te^2 v} dx\right)^{1/2}\left(\int_0^1 1_{(\te\le \bar\te)} {\te v} dx\right)^{1/2}\\ &\le C V^{1/2}(t)   ,\ea\enn
which implies that for   $t>0,$ \be \la{twq1}\min_{x\in[0,1]}\te(x,t)\ge \frac{\alpha_1}{4}-CV(t ).\ee
Combining  \eqref{cq8}--\eqref{twq1}  yields that
\be\ba\la{cq12a} v(x,t)&\ge C^{-1} \int_0^t e^{-2(t-\tau)}\min_{x\in[0,1]}\te(x,\tau)d\tau \\&\ge  C^{-1}\int_0^t e^{-2(t-\tau)}\left(\frac{\al_1}{4}-C  V( \tau)\right)d\tau\\ &\ge \frac{ C^{-1} \al_1}{8} - \frac{ C^{-1}\al_1}{8} e^{-2t} -C   \int_0^t e^{-2(t-\tau)} V( \tau) d\tau.\ea\ee
Since
\bnn\ba  \int_0^t e^{-2(t-\tau)} V( \tau) d\tau&=\int_0^{t/2} e^{-2(t-\tau)} V( \tau) d\tau+\int_{t/2}^t e^{-2(t-\tau)} V( \tau) d\tau\\&\le e^{-t}\int_0^\infty V(\tau)d\tau +\int_{t/2}^t V( \tau) d\tau\rightarrow 0, \mbox{ as }t\rightarrow \infty,\ea\enn it follows from \eqref{cq12a} that  there exists some $\ti T >0$ such that
\be\ba\la{cq12} v(x,t) \ge  \frac{ C^{-1} \al_1}{16} \ea\ee for all $(x,t)\in[0,1]\times[ \ti T ,+\infty).$

Finally,    using \eqref{2.6}, \eqref{cq8}, and \eqref{cq7},  we get $$ v(x,t) \ge   C^{-1}e^{-2\ti T},$$ for all $(x,t)\in[0,1]\times[0, \ti T ],$  which together with \eqref{cq12} implies \bnn\la{cq122} v(x,t) \ge C_0\triangleq \min\left\{\frac{ C^{-1} \al_1}{16},   C^{-1}e^{-2\ti T}\right\},\enn  for all $(x,t)\in[0,1]\times[ 0 ,+\infty).$ We finish the proof of Lemma \ref{wq20}. \thatsall

To obtain the upper bound of $v,$  we
set \bnn M_v(t)\triangleq 1+\max_{x\in[0,1]}v(x,t).\enn Then we have the following time-independent upper bound of $v$ for $\b>0.$

\begin{lemma} \la{lemq2} For $\beta>0,$ there exists a   positive constant $C $  such that for all   $(x,t)\in[0,1]\times [0,+\infty),$  \be \la{vg14}v(x,t)\le C.\ee\end{lemma}

\pf First, we claim that
\be \la{u1}\int_0^T\max_{x\in[0,1]}\left(\te^{\frac{ 1}{2}}(x, t)- 2^{\frac{ 1}{2}}  \right)_+^2 dt \le C     . \ee
Indeed, on the one hand, for $\beta\in [1,\infty),$ we have
\be\la{u2}\ba   &\int_0^T\max_{x\in[0,1]}\left(\te^{\frac{ 1}{2}}(x, t)- 2^{\frac{ 1}{2}}  \right)_+^2 dt\\ &\le C\int_0^T\max_{x\in[0,1]}\left(\te^{\frac{\b}{2}}(x, t)- 2^{\frac{\b}{2}}  \right)_+^2dt \\&\le C\int_0^T\left(\int_0^1 \left|\pa_x\left(\te^{\frac{\b}{2}}(x, t)- 2^{\frac{\b}{2}}  \right)_+\right|dx\right)^2dt\\&\le C\int_0^T \int_0^1 \frac{ \te^\b \te_x^2}{\te^2 v} dx \int_0^1  v dxdt  \\ &\le C    .\ea\ee
On the other hand, for $\beta\in (0,1) $ and   $\eta\triangleq  (2-\b)/4 \in\left(0,\frac{1}{2}\right)$, integrating \eqref{2.4} multiplied by
 $\left(\te^{\eta}-4^{\eta}\right)_{+}\te^{\eta-1}$ over $(0,1)\times(0,T)$, we get
 \be\ba\la{2.24}&\frac{\beta}{2}\int_0^T\int_{(\te>4)(t)} \frac{\te^{-1+\b/2}\te_x^2}{v }dxdt\\&\quad+\int_0^T\int_{\O}\frac{  u_x^2+|\nw_x|^2+|\nb_x|^2}{v}
 \left(\te^\eta-4^\eta\right)_{+}\te^{\eta-1} dxdt\\&
 =\frac{1}{2\eta}\int_{\O}\left(\left(\te^\eta-4^\eta\right)_{+}^2-\left(\te_0^\eta-4^\eta\right)_{+}^2\right)dx+4^\eta(1-\eta)\int_0^T\int_{(\te>4)(t)}\frac{\te^\b\te_x^2}{v\te^{2-\eta}}dxdt\\&\quad
 +\int_0^T\int_{\O}\frac{\te u_x}{v}\left(\te^\eta-4^\eta\right)_{+}\te^{\eta-1}dxdt\\&
 \le C(\ve)+  \ve\int_0^T\int_{(\te>4)(t)} \frac{\te^{-1+\b/2}\te_x^2}{v }dxdt+\frac{1}{2}\int_0^T\int_{\O}\frac{  u_x^2}{v}\left(\te^\eta-4^\eta\right)_{+}\te^{\eta-1}dxdt  ,   \ea\ee
where in the last inequality we have used
\bnn\ba &\int_0^T\int_0^1\frac{1}{  v}\left(\te^\eta-4^\eta\right)_{+}\te^{\eta+1} dxdt\\&
\le C\int_0^T\max_{x\in [0,1]}(\te^{1/2}(x,t)-2^{1/2})^2_+\int_0^1 \te^{2\eta }dxdt\\&
\le  \varepsilon\int_0^T\int_{(\te>4)(t)} \frac{\te^{-1+\b/2}\te_x^2}{v }dxdt+C(\varepsilon), \ea\enn
due to
\be \la{2.26}\ba &\int_0^T\max_{x\in [0,1]}(\te^{1/2}(x,t)-2^{1/2})^2_+ dt  \\&
\le  C\int_0^T\int_0^1 \frac{\te_x^2}{\te v }dx\int_0^1 v  dxdt \\&
\le \varepsilon\int_0^T\int_0^1 \frac{\te^{-1+\b/2}\te_x^2}{v }dxdt+C(\varepsilon) \int_0^T\int_0^1 \frac{\te^{-2 +\b }\te_x^2}{v }dxdt   \\&
\le \varepsilon\int_0^T\int_{(\te>4)(t)} \frac{\te^{-1+\b/2}\te_x^2}{v }dxdt+C(\varepsilon). \ea\ee
Combining   \eqref{2.24} with \eqref{2.26}  gives \eqref{u1} for $\b\in (0,1)$.

Next, it follows from \eqref{wq1}  that
\bnn\ba\la{bqq2}     \left|\te^{\frac{ 1}{2}}(x,t)- \bar\te^{\frac{ 1}{2}}( t) \right|&\le\left|\te^{\frac{\b+1}{2}}(x,t)- \bar\te^{\frac{\b+1}{2}}( t) \right|\\  &\le \frac{\b+1}{2}\left(\int_0^1 \frac{ \te^\b \te_x^2}{\te^2 v} dx\right)^{1/2}\left(\int_0^1  {\te v} dx\right)^{1/2}\\ &\le C V^{1/2}(t)  M_v^{1/2}(t),\ea\enn which together with \eqref{qwe1} leads to
\be\ba\la{bqq7}
   \theta(x,t)\le C +CV(t)M_v(t),
\ea\ee
for all $(x,t)\in [0,1]\times [0,\infty).$

Finally, standard calculations give \be \la{44}\ba   \max_{x\in[0,1]}|\nb|^2(x,t) &\leq C \int_0^1|\nb\cdot\nb_x| dx  \\&\leq C \int_0^1\frac{|\nb_x|^2}{v\te}dx +C \int_0^1v\te|\nb|^2dx  \\&\leq CV(t) +C \max_{x\in[0,1]}(\te^{1/2} (x,t)-2^{1/2})_+^2+C  ,\ea\ee which together with  \eqref{2.6}, \eqref{cq8},  \eqref{cq7},  and \eqref{bqq7} leads to
\bnn\ba\la{cq10} v(x,t)
&\le C+C\int_0^t e^{-\al_1(t-\tau)}\max_{x\in[0,1]} \left(\te +v|\nb|^2\right)(x,\tau)d\tau\\
&\le C+C\int_0^t e^{-\al_1(t-\tau)}\left(1+(1+ V(\tau))M_v(\tau) \right)d\tau\\&\quad +C\int_0^t  \max_{x\in[0,1]}(\te^{1/2} (x,\tau)-2^{1/2})_+^2 M_v(\tau) d\tau .\ea\enn
We thus  obtain   \eqref{vg14} from this, \eqref{2.2}, \eqref{u1}, and the Gronwall inequality. The proof of Lemma \ref{lemq2} is finished.
\thatsall

For $\b>0,$ to obtain the  time-independent lower bound of the temperature, we  need the following uniform (with respect to time) estimate on the $L^\infty(0,T;L^p)$-norm of $\te^{-1}.$

\begin{lemma} \la{len1} For $\beta> 0 $ and any $p>0,$   there exists some positive constant $C(p)$  such that \be \la{ljj1}   \sup_{\le t\le T}\int_0^1\te^{-p}dx+ \int_0^T\int_0^1   \te^{\b-1-p} \te_x^2    dxdt  \le C(p) .\ee \end{lemma}

 \pf
It   suffices to prove \eqref{ljj1} for $p\not=1$ since it holds for $p=1 $ due to \eqref{2.2}.
Multiplying \eqref{2.4}   by $1/\te^p$    and integration by parts gives
\be\ba\la{aa1}  &\frac{1}{p-1}\left(\int_0^1\te^{1-p}dx\right)_t+p\int_0^1 \frac{\te^\b\te_x^2}{v\te^{p+1}}dx+\int_0^1\frac{u_x^2+|\nw_x|^2+|\nb_x|^2}{v\te^p}dx \\& =\int_0^1 \frac{\left(\te^{1-p}-1\right)u_x}{v}dx+\int_0^1 \frac{ u_x}{v}dx\\ & \le C(p)\int_0^1   \left|\te^{\frac{1}{2}}-1 \right|\left(\te^{\frac{1}{2}-p}+1\right)| u_x| dx +\left(\int_0^1 \ln v dx\right)_t  \\ & \le  C(p)\max_{x\in[0,1]}\left| \te^{  \frac{1}{2}}-1\right| \left( \int_0^1 \te^{1-p}dx \right)^{1/2}\left( \int_0^1\frac{u_x^2}{v\te^p}dx \right)^{1/2}   \\&\quad  +\max_{x\in[0,1]}\left| \te^{  \frac{1}{2}}-1\right|\int_0^1 |u_x|dx+\left(\int_0^1 \ln v dx\right)_t\\ & \le  C(p)\max_{x\in[0,1]}\left( \te^{  \frac{1}{2}}-1\right)^2  \left(1+ \int_0^1 \te^{1-p}dx \right)  + \frac12 \int_0^1\frac{u_x^2}{v\te^p}dx \\&\quad  + C(p)\left(\int_0^1 |u_x|dx\right)^2 +\left(\int_0^1 \ln v dx\right)_t.\ea\ee

 Next, it follows from \eqref{wq1}  and \eqref{qwe1} that $$\alpha_1\le\int_0^1 \te dx\le \int_0^1\left(\te+\eta\frac{u^2+|\nw|^2+v|\nb|^2}{2}\right)dx\le 1,$$ which yields that  for any real number $q,$
\be\ba\la{ux} | 1-\bar \te^{q}| &=\left|\int_0^1 \frac{d}{d\eta } \left(\int_0^1\left(\te+\eta  \frac{u^2+|\nw|^2+v|\nb|^2}{2}\right)dx\right)^{q}   d\eta\right|  \\&\le C(q) \int_0^1 \left({u^2} +|\nw|^2+v|\nb|^2\right) dx \\& \le  C(q)\max_{x\in[0,1]} \left(|u|+|\nw| + |\nb|\right)\left(\int_0^1 \left(u^2+|\nw|^2+v|\nb|^2\right)   dx\right)^{1/2} \\& \le C \int_0^1 \left(|u_x|+|\nw_x| + |\nb_x|\right) dx  \\ &\le  C\left( \int_0^1\frac{ u_x^2+ |\nw_x|^2 + |\nb_x|^2}{v\te}dx\right)^{1/2}\left(\int_0^1 v\te dx\right)^{1/2}\\ &\le   CV^{1/2}(t). \ea\ee

Then,   it follows from \eqref{ux} and \eqref{qwe1} that for $ \b\in (0,1),$
\be\la{qwq2}\ba \max_{x\in[0,1]}\left|\te^{\frac{1 }{2}}-1\right|&\le \max_{x\in[0,1]}\left|\te^{\frac{1 }{2}}-\bar\te^{\frac{1 }{2}}\right| +\max_{x\in[0,1]}\left|\bar\te^{\frac{1 }{2}}-1\right| \\ &\le  C  \int_0^1  \te^{ -\frac{1 }{2} }|\te_x|dx+C V^{1/2}(t)\\&\le C\left(\int_0^1  \te^{\b-2}\te_x^2dx\right)^{1/2}\left(\int_0^1  \te^{1-\b} dx\right)^{1/2}+ CV^{1/2}(t)\\&\le CV^{1/2}(t),\ea\ee
and that for $\b\ge 1,$
 \be\la{qw2}\ba \max_{x\in[0,1]}\left|\te^{\frac{1 }{2}}-1\right|&\le \max_{x\in[0,1]}\left|\te^{\frac{1 }{2}}-\bar\te^{\frac{1 }{2}}\right| +\max_{x\in[0,1]}\left|\bar\te^{\frac{1 }{2}}-1\right|\\ &\le C \max_{x\in[0,1]}\left|\te^{\frac{\b }{2}}-\bar\te^{\frac{\b}{2}}\right| +CV^{1/2}(t) \\ &\le  C  \int_0^1  \te^{ \frac{\b }{2}-1}|\te_x|dx+C V^{1/2}(t)\\&\le  CV^{1/2}(t). \ea\ee
It thus follows from \eqref{qwq2}, \eqref{qw2}, and \eqref{2.2} that for $\b> 0,$
 \be\la{qw2s}\ba\int_0^T \max_{x\in[0,1]}\left(\te^{\frac{1 }{2}}-1\right)^2dt\le C.\ea\ee

Finally, noticing  that  for $p\in [0,1],$
\bnn\la{qwq1}\ba
\int_0^1\te^{1-p}dx\le \int_0^1\te dx+ 1\le C,
\ea\enn and that both
  \eqref{2.2} and   \eqref{wq1} imply that for $\b\ge 0,$
\be\la{22eq5}
\sup_{0\le t<\infty}\int_0^1  |\ln v|dx\le C,
\ee
 after using \eqref{ux}, \eqref{2.2}, \eqref{qw2s},   and the Gronwall  inequality,   we   obtain  \eqref{ljj1} from \eqref{aa1} and   finish  the proof of Lemma \ref{len1}.
\thatsall

\begin{lemma} \la{lemq1} For $\beta\ge 0,$ there exists a   positive constant $C $  such that for all $T> 0,$ \be \la{uiz1}\sup_{0\le t\le T}\int_0^1v^2_xdx +\int_0^T\int_0^1\left( (1+\te) v_x^2+u_x^2+|\nb_x|^2+|\nw_x|^2 \right)dxdt\le C.\ee \end{lemma}

\pf
{\it Case 1 $(\b=0).$} First,
%\be \int_0^1(v\nb)(x,t)dx\equiv \ti\nb_0=(\ti b^1_0,\ti b^2_0)=  \int_0^1 v_0\nb_0 dx,\ee   which shows that for $t>0$ there exists some points $x_i(t) (i=1,2)$ such that \be (b^i-\ti b^i_0)(x_i(t),t)=0, i=1,2\ee
  multiplying \eqref{1.3} by $\nw$ and integrating the resulting equality over $(0,1)  $ yields
\be\la{vd1}\ba &\frac{1}{2}\frac{d}{dt} \int_0^1 |\nw|^2dx +\int_0^1 \frac{|\nw_x|^2}{v}dx%\\&=-\int_0^1 \left(\nb-\ti \nb_0 \right) \cdot \nw_x dx
\\&\le C\int_0^1 |\nw_x|\left|\nb  \right|dx\\& \le C\int_0^1|\nw_x| | \nb   |  (\ln \te-\ln \bar\te)_+dx+C\int_0^1 |\nw_x| | \nb |(1-\ln \te)_+dx\\& \le C\left(\int_0^1 \frac{|\nw_x|^2}{v}dx\right)^{1/2} \left(\int_0^1 v |\nb |^2  dx\right)^{1/2}\max_{x\in[0,1]} (\ln \te-\ln \bar\te)_+\\&\quad +C\left(\int_0^1 \frac{|\nw_x|^2}{v\te }dx\right)^{1/2} \max_{x\in[0,1]} (| \nb |  (\ln 2-\ln \te)_+)  \\& \le \frac12 \int_0^1 \frac{|\nw_x|^2}{v}dx+CV(t),
\ea\ee
where in the last inequality we have used the following two simple facts:
\bnn\ba \max_{x\in[0,1]} (\ln \te-\ln \bar\te)_+&\le C \int_0^1|((\ln \te-\ln \bar\te)_+)_x|dx\\&\le C \left(\int_0^1v dx \int_0^1  \frac{\te_x^2}{v \te^2}dx\right)^{1/2} \\&\le CV^{1/2}(t),\ea\enn
and
\bnn\ba &\max_{x\in[0,1]} (| \nb|  (1-\ln \te)_+)\\ &\le C\int_0^1 |  (\nb      (1-\ln \te)_+)_x|dx\\&\le C\int_0^1 | \nb_x | \te^{-1/2}dx+  C\int_0^1 |  \nb   |  \frac{| \te_x|}{\te}dx\\&\le C\left(\int_0^1  \frac{ \nb_x^2}{v \te}dx \int_0^1vdx\right)^{1/2}+ C\left(\int_0^1v|\nb |^2dx \int_0^1  \frac{\te_x^2}{v \te^2}dx\right)^{1/2}\\&\le CV^{1/2}(t).\ea\enn Denoting \bnn \ti V(t)\triangleq \int_0^1 \frac{|\nw_x|^2}{v}dx+V(t),\enn
We thus obtain from \eqref{vd1} and \eqref{2.2} that
 \be\la{q2.45}\int_0^T\ti V(t) dt  \le C .\ee

Next, using \eqref{1.1}, we rewrite \eqref{1.2} as \be\la{k1}  (\ln v)_{xt}=u_t+\left(\frac{\te}{v}\right)_x+\nb\cdot \nb_x.\ee
Adding  \eqref{k1}  multiplied by $(\ln v)_x$ to  \eqref{1.4}    by $v\nb,$   and integrating the resulting equality over $(0,1)\times(0,T),$
 one has \be\la{vg1}\ba &\frac12 \frac{d}{dt}\int_0^1 \left((\ln v)_x^2+v^2|\nb|^2\right)dx+\int_0^1 \left(  |\nb_x|^2+\frac{\te}{v}(\ln v)_x^2\right)dx\\& =\frac{d}{dt} \int_0^1 u(\ln v)_xdx +\int_0^1\frac{u_x^2}{v}dx+\int_0^1 \frac{\te_x(\ln v)_x }{v}dx+\int_0^1 v\nw_x\cdot\nb dx.\ea\ee

Then, on the one hand,
 for the fourth  term on the righthand side of \eqref{vg1}, we have by \eqref{vd1},
\be\la{vg12}\ba   \left|\int_0^1 v\nw_x \cdot \nb dx\right| & \le CM_v(t)\int_0^1   |\nw_x|| \nb | dx  \\& \le  CM_v(t)\ti V(t).\ea\ee

On the other hand, for the third term on the righthand side of \eqref{vg1}, integrating by part  gives
\be\la{vg100}\ba   \int_0^1 \frac{\te_x(\ln v)_x }{v}dx   &=-\int_0^1  \ln \frac{v}{C_0} \left(\frac{\te_x }{v}\right)_xdx\\&=-\int_0^1  \ln \frac{v}{C_0} \left(
\theta_{t}+  \frac\theta v
u_{x}-\frac{  u_{x}^{2}+|\nw_x|^2+|\nb_x|^2}{v}\right) dx\\&=-\left(\int_0^1 \theta \ln \frac{v}{C_0} dx\right)_{t} +\int_0^1 \frac{\theta u_{x}}{ v
} dx\\&\quad -\int_0^1 \frac{\theta u_{x}}{ v
}   \ln \frac{v}{C_0}  dx+\int_0^1 \frac{  u_{x}^{2}+|\nw_x|^2+|\nb_x|^2}{v}\ln \frac{v}{C_0} dx . \ea\ee
For the second term on the righthand side of \eqref{vg100}, we have
\be \la{vg3} \ba    \int_0^1 \frac{\te }{v}u_xdx   &=\int_0^1 \frac{\left(\te-1\right)}{v}u_xdx + \int_0^1 \frac{ u_x}{v}dx \\&\le \ve \int_0^1 \frac{u_x^2}{v}dx+C(\ve)\int_0^1 \left(\te-1\right)^2dx    + \left(\int_0^1 \ln vdx\right)_t \\&\le \ve\int_0^1 \frac{u_x^2}{v}dx +  C(\ve) V(t)M_v(t)  + \left(\int_0^1 \ln vdx\right)_t,\ea\ee
where in the last inequality we have used
\be\la{q2.32}\ba  \int_0^1 \left(\te-1\right)^2dx   &\le C\max_{x\in[0,1]} \left(\te^{1/2}-1\right)^2  \\&\le C\max_{x\in[0,1]} \left(\te^{1/2}-\bar\te^{1/2}\right)^2+C(1-\bar\te)^2 \\&\le C \left( \int_0^1|\te_x|\te^{-1/2}dx \right)^2+CV(t)M_v(t) \\&\le C  \int_0^1 \frac{\te_x^2}{ \te^{2}v}dx  \int_0^1 v\te dx+CV(t)M_v(t) \\& \le CV(t)M_v(t)\ea\ee
due to
\be\la{vg4}\ba   1-\bar \te  &=\frac12   \int_0^1   (u^2+|\nw|^2+v|\nb|^2) dx        \\&\le C \int_0^1 \left(|u_x|+|\nw_x| + |\nb_x|\right) dx\\  &\le  C\left( \int_0^1\frac{ u_x^2+|\nw_x|^2+|\nb_x|^2}{v\te}dx\right)^{1/2}\left(\int_0^1 v\te dx\right)^{1/2}\\ &\le   CV^{1/2}(t)M^{1/2}_v(t). \ea\ee
Similarly,
for the third term on the righthand side of \eqref{vg100}, we have
\be \la{vg11} \ba  &  -\int_0^1 \frac{\te }{v}u_x\ln \frac{v}{C_0}dx  \\&\le C\left(\int_0^1 \frac{u_x^2}{v}dx  +     V(t)M_v(t)\right) \ln M_v(t )  -\frac12 \left(\int_0^1\ln^2 \frac{v}{C_0}dx\right)_t,\ea\ee
due to \eqref{vg3}. Putting \eqref{vg3} and \eqref{vg11} into \eqref{vg100} gives
\be\la{vg10}\ba   \int_0^1 \frac{\te_x(\ln v)_x }{v}dx    &\le \frac{d}{dt} \int_0^1\left(- \theta \ln \frac{v}{C_0} + \ln v-\frac12\ln^2 \frac{v}{C_0}\right) dx  \\&\quad  + C \left(\int_0^1 \frac{u_x^2+|\nb_x|^2}{v}dx   + \ti     V(t) M_v(t)\right)\ln M_v(t). \ea\ee

Putting   \eqref{vg12}  and \eqref{vg10} into \eqref{vg1} gives
  \be\la{vg13}\ba &\frac12 \frac{d}{dt}\int_0^1 \left((\ln v)_x^2+|v\nb|^2+2\te \ln\frac{v}{C_0}+\ln^2\frac{v}{C_0} -2\ln v-2u(\ln v)_x\right)dx\\&+\int_0^1 \left(  |\nb_x|^2+\frac{\te}{v}(\ln v)_x^2\right)dx \\&\le    C \left(\int_0^1 \frac{u_x^2+|\nb_x|^2}{v}dx   +      \ti  V(t) M_v(t)\right)\ln M_v(t).\ea\ee

Then, to estimate the first term on the righthand side of \eqref{vg13}, we multiply \eqref{1.2} by $u$ and integrate the resultant equality over $(0,1)$ to get
\be  \la{vg2}\ba   \frac12\left(\int_0^1 u^2dx\right)_t+\int_0^1 \frac{u_x^2}{v}dx  =\int_0^1 \frac{\te }{v}u_xdx-\int_0^1 \nb \cdot \nb_x  u dx . \ea\ee

Thus, it follows from \eqref{vg4} and \eqref{2.2} that
\be\la{vg6}\ba \left|\int_0^1  \nb\cdot \nb_x  u dx\right|&\le \max_{x\in[0,1]}|u|\left(\int_0^1 v|\nb|^2\  dx\right)^{1/2} \left(\int_0^1 \frac{|\nb_x|^2}{v } dx\right)^{1/2} \\&\le \ve\int_0^1 \frac{|\nb_x|^2}{v}dx +C(\ve)V(t)M_v(t),\ea\ee
which together with \eqref{vg2} and \eqref{vg3} gives
\be \la{vg5} \ba & \frac12 \left(\int_0^1 u^2dx\right)_t+\frac12 \int_0^1 \frac{u_x^2}{v}dx   \\&\le  \frac14\int_0^1 \frac{|\nb_x|^2}{v}dx+  C V(t)M_v(t)  + \left(\int_0^1 \ln vdx\right)_t.\ea\ee
Multiplying \eqref{1.4} by $\nb$ and integrating the result over $(0,1),$ we obtain from \eqref{vg6} that
\bnn  \la{vg7}\ba & \frac12\left(\int_0^1 v|\nb|^2dx\right)_t+\int_0^1 \frac{|\nb_x|^2}{v}dx \\&=\int_0^1 \nb\cdot\nb_xu dx-\int_0^1 \nw\cdot\nb_xdx\\&\le \frac14\int_0^1 \frac{|\nb_x|^2}{v}dx +C V(t)M_v(t)+C \max_{x\in [0,1]}|\nw|^2 \\&\le \frac14\int_0^1 \frac{|\nb_x|^2}{v}dx +C V(t)M_v(t)+C\int_0^1 \frac{|\nw_x|^2}{v}dx , \ea\enn which combined with
  \eqref{vg5} gives
\bnn  \la{vg8}\ba & \left(\int_0^1 (u^2+v|\nb|^2)dx\right)_t+\int_0^1 \frac{u_x^2+|\nb_x|^2}{v}dx   \\&\le    C \ti V(t)M_v(t)  + 2\left(\int_0^1 \ln vdx\right)_t.\ea\enn
Integrating this over $(0,t),$ we have by \eqref{22eq5}
\be \la{vg9}\int_0^t \int_0^1 \frac{u_x^2+|\nb_x|^2}{v}dx  ds\le C+C\int_0^t\ti V(s)M_v (s)ds.\ee

Next, it follows from  \eqref{vg13} and \eqref{vg9} that
\be\la{2.45}\ba & \sup_{0\le s\le t}\int_0^1  (\ln v)_x^2dx +\int_0^t \int_0^1 \left(  |\nb_x|^2+\frac{\te}{v}(\ln v)_x^2\right)dxds \\&\le  C\ln \sup_{0\le s\le t}M_v(s) + C \int_0^t \ti V(s) M_v(s)ds \ln \sup_{0\le s\le t}M_v(s) \\&\le  C(\ve)\left(1+ \int_0^t \ti V(s) M_v(s)ds  \right)\ln \left(1+ \int_0^t \ti V(s) M_v(s)ds  \right) \\&\quad+  C(\ve)+\ve   \sup_{0\le s\le t}M_v(s),\ea\ee where in the second inequality  we have used  \bnn fg\le  e^f-f-1+(1+g)\ln(1+g)-g, \mbox{ for any } f\ge 0,g\ge 0, \enn with $$f=\frac12 \ln \sup_{0\le s\le t}M_v(s), \quad g= 2C \int_0^t \ti V(s) M_v(s)ds.$$

Then, direct computation shows
\bnn\ba  v-1 &\le C\left(\int_0^1v^2 dx\right)^{1/2}\left(\int_0^1\frac{v_x^2}{v^2}dx\right)^{1/2} \\ &\le C M^{1/2}_v(t) \left(\int_0^1\frac{v_x^2}{v^2}dx\right)^{1/2}, \ea\enn which gives \bnn  M_v(t)\le C+C\int_0^1 (\ln v)_x^2 dx.\enn
Combining this, \eqref{2.45},  \eqref{q2.45}, and the Gronwall inequality shows that
 for any $(x,t)\in[0,1]\times [0,+\infty),$
\be \la{q2.49}v(x,t)\le C,\ee which together with \eqref{2.45}   and \eqref{q2.45} implies \be \la{uiz11}\sup_{0\le t\le T}\int_0^1v^2_xdx +\int_0^T\int_0^1\left(   \te  v_x^2+u_x^2+|\nb_x|^2+|\nw_x|^2 \right)dxdt\le C.\ee

Finally, direct computation shows
\be\la{pz1}\ba    \int_0^1 v_x^2dx&= \int_0^1 v_x^2 \left(1-\te\right) dx+  \int_0^1 v_x^2\te dx\\&  \le \int_0^1 v_x^2 \left(1-\te\right)_+ dx+  \int_0^1 v_x^2\te dx \\& \le \frac{1}{2}\int_0^1 v_x^2  dx +C\max_{x\in[0,1]}\left(\te^{1/2}-1\right)^2 +  \int_0^1 v_x^2\te dx,\ea\ee where in the last inequality we have used \eqref{uiz11}. Combining this, \eqref{uiz11}, \eqref{q2.32}, and \eqref{q2.49} gives \eqref{uiz1} for $\b=0$.

{\it Case 2 $(\b>0).$}   First, we rewrite the momentum equation \eqref{1.2} as
\bnn\la{bbbb}\ba \left(u-\frac{  v_x}{v}\right)_t=-\left(\frac{\te}{v}+\frac{1}{2}|\nb|^2\right)_x.   \ea\enn
Multiplying the above equation by $u-\frac{  v_x}{v}$ and integrating the resultant equality  yields that for any $t\in(0,T)$
\be\ba\la{cccc}&\frac{1}{2}\int_0^1\left(u-\frac{  v_x}{v}\right)^2 dx -\frac{1}{2}\int_0^1\left(u-\frac{  v_x}{v}\right)(x,0)dx\\&= \int_0^t \int_0^1\left( \frac{\te v_x}{v^2}-\frac{\te_x}{v}-\nb\cdot\nb_x\right)\left(u-\frac{  v_x}{v}\right)dxdt \\
&=  -\int_0^t\int_0^1\frac{  \te v_x^2}{v^3}dxdt+\int_0^t\int_0^1\frac{\te u v_x}{v^2}dxdt \\&\quad-\int_0^t\int_0^1\frac{\te_x}{v} \left(u-\frac{  v_x}{v}\right)dxdt-\int_0^t\int_0^1\nb\cdot\nb_x  \left(u-\frac{  v_x}{v}\right)dxdt \\
&=  -\int_0^t\int_0^1\frac{  \te v_x^2}{v^3}dxdt+\sum_{i=1}^3I_i.\ea\ee
Each $I_i (i=1,2,3)$ can be estimated as follows:

First, Cauchy's inequality gives
\be\ba\la{dddd}|I_1|&\leq\frac{1}{4}\int_0^t\int_0^1\frac{  \te v_x^2}{v^3}dxdt+C\int_0^T\int_0^1 {u^2\te} dxdt   \\&\leq \frac{1}{4}\int_0^t\int_0^1\frac{  \te v_x^2}{v^3}dxdt+C,\ea\ee
where  we have used
\be\ba\la{ddad}   \int_0^T\int_0^1 {u^2\te} dxdt  &\leq C\int_0^T\int_0^1 u^2((\te^{1/2}-1)^2+1)dxdt \\&\leq C\int_0^T\left(\int_0^1u^2dx+\max_{x\in[0,1]}\left(\te^{1/2}-1\right)^2\right) dt \\&\leq C ,\ea\ee due to  \eqref{ux}, \eqref{2.2}, and \eqref{qw2s}.

Next, using \eqref{ddad}, \eqref{2.2}, and \eqref{ljj1} with $p=\b$, we have
\be\ba\la{eeee}|I_2| &\leq C\int_0^T\int_0^1 {u^2\te} dxdt+C\int_0^T\int_0^1 \te^{-1}\te_x^2 dxdt+\frac{1}{2}\int_0^t\int_0^1\frac{  \te v_x^2}{v^3}dxdt \\
&\leq C+ \frac{1}{2}\int_0^t\int_0^1\frac{  \te v_x^2}{v^3}dxdt.\ea\ee

 Finally, integrating \eqref{2.4} over $(0,1)\times (0,T),$ we have by \eqref{22eq5}
\bnn \ba&\int_0^T\int_0^1\frac{  u_x^2+|\nw_x|^2+|\nb_x|^2}{v}dxdt\\ &= \int_0^1\te dx-\int_0^1\te_0dx +\int_0^T\int_0^1\frac{\te-1}{v}u_xdx+\int_0^1\ln vdx-\int_0^1\ln v_0dx\\&\le C+\frac{1}{2}\int_0^T\int_0^1\frac{  u_x^2 }{v}dxdt +C \int_0^T\max_{x\in[0,1]}(\te^{1/2}-1)^2 dt , \ea\enn
which together with \eqref{qw2s}  gives
\be\la{h}\ba&\int_0^T\int_0^1(u_x^2+|\nw_x|^2+|\nb_x|^2)dxdt  \le C . \ea\ee Combining this with Cauchy's inequality leads to
\be\ba\la{gggg}|I_3|&\leq C\int_0^t\int_0^1\left(|\nb_x|^2+|\nb|^2 \left(u-\frac{  v_x}{v}\right)^2\right)dxdt \\&\leq C+C\int_0^tV(t)\int_0^1\left(u-\frac{  v_x}{v}\right)^2dxdt,\ea\ee due to \bnn \max_{0\le x\le 1}|\nb|^2\le C\int_0^1\frac{|\nb_x|^2}{v\te}dx\int_0^1v\te dx\le CV(t).\enn
Putting \eqref{dddd}, \eqref{eeee}, and \eqref{gggg} into \eqref{cccc}, we  obtain after using the Gronwall  inequality and  \eqref{2.2}  that
\bnn\la{nnnn}\ba \sup_{0\le t\le T}\int_0^1\left(u-\frac{  v_x}{v}\right)^2dx+\int_0^T\int_0^1\frac{\te v_x^2}{v^3}dxdt\leq C,\ea\enn
which together with \eqref{2.2}  gives
 \be\la{vvx1}\ba
\sup_{0\le t \le T}\int_0^1 v_x^2  dx+\int_0^T\int_0^1 v_x^2\te dxdt\le C.
\ea\ee

Finally,  since \eqref{pz1} still holds for $\b>0$ due to \eqref{vvx1},  we have
\bnn \int_0^T\int_0^1 v_x^2  dxdt\le C,\enn which  together with  \eqref{vvx1} and \eqref{h} proves \eqref{uiz1}   for  $\b>0.$  The proof of Lemma \ref{lemq1} is  finished. \thatsall

\begin{lemma}\la{lemmy0} For $\b\ge 0,$
there is a positive constant $C$ such that for all $T> 0,$
\be \la{aas1}\ba
&\sup_{0\le t\le T} \int_0^1\left( |\nb_x|^2+|\nw_x|^2\right)dx \\&  +\int_0^T\int_0^1\left(  |\nb_t|^2 +|\nb_{xx}|^2 +|\nw_t|^2+|\nw_{xx}|^2\right)dx dt\leq C.\ea\ee
\end{lemma}

\pf First,   rewriting  \eqref{1.3} as
\be\la{ppp}\ba \nw_t=\frac{\nw_{xx}}{v}-\frac{\nw_xv_x}{v^2}+\nb_x, \ea\ee multiplying \eqref{ppp} by $\nw_{xx},$ and integrating the resulting equality over $(0,1)\times(0,T)$, we obtain after using \eqref{1.9},   \eqref{uiz1}, and Cauchy's inequality that
\be\ba\la{bbb}&\frac{1}{2}\int_0^1|\nw_x|^2dx+\int_0^T\int_0^1\frac{|\nw_{xx}|^2}{v}dxdt   \\& \leq C+\frac{1}{2}\int_0^T\int_0^1\frac{|\nw_{xx}|^2}{v}dxdt +C\int_0^T\int_0^1\left(|\nb_x|^2+|\nw_x|^2v_x^2\right)dxdt \\&\leq C+\frac{1}{2}\int_0^T\int_0^1\frac{|\nw_{xx}|^2}{v}dxdt +C\int_0^T\max_{x\in[0,1]}|\nw_x|^2dt.\ea\ee

Noticing that    for any $f\in \{   \int_0^1fdx=0\}\cup \{  f(0)=0\},$
\be\ba\la{ux2}
\max_{x\in[0,1]}f^2 \le 2\left(\int_0^1 f^2d x\right)^{1/2}\left(\int_0^1 f_x^2dx\right)^{1/2},
\ea\ee  we get for any $\ve>0,$
\be\ba\la{ccc}\int_0^T\max_{x\in[0,1]}|\nw_x|^2dt&\leq C(\ve)\int_0^T\int_0^1|\nw_x|^2dxdt+\varepsilon\int_0^T\int_0^1\frac{|\nw_{xx}|^2}{v}dxdt \\&\leq C(\ve)+\varepsilon\int_0^T\int_0^1\frac{|\nw_{xx}|^2}{v}dxdt,\ea\ee
which combined with \eqref{bbb}  leads to
\bn\ba \la{ddd}\sup_{0\le t\le T}\int_0^1|\nw_x|^2dx+\int_0^T\int_0^1|\nw_{xx}|^2dxdt\leq C.\ea\en
Combining this, \eqref{ppp},  \eqref{ccc}, and \eqref{uiz1}  gives
\be\ba\la{ooo}\int_0^T\int_0^1|\nw_t|^2dxdt &\leq
C\int_0^T\int_0^1 \left(|\nb_x|^2+|\nw_{xx}|^2 + v_x^2|\nw_x|^2\right)dxdt\\&\leq  C+ C\int_0^T\max_{x\in[0,1]}|\nw_x|^2dt\\&\leq C .\ea\ee

Next,   rewriting \eqref{1.4} as
\be\la{kkk}\ba\nb_t=\frac{\nw_x}{v}+\frac{\nb_{xx}}{v^2}-\frac{\nb_x v_x}{v^3}-\frac{\nb u_x}{v},\ea\ee multiplying \eqref{kkk} by $\nb_{xx}$ and integrating the result over $(0.1)\times(0,T)$, we deduce from \eqref{uiz1},  \eqref{2.2}, \eqref{ux2},      and Cauchy's inequality that
\bnn\ba &\frac{1}{2}\int_0^1|\nb_x|^2dx+\int_0^T\int_0^1\frac{|\nb_{xx}|^2}{v^2}dxdt\\ &\leq C+\frac{1}{2}\int_0^T\int_0^1\frac{|\nb_{xx}|^2}{v^2}dxdt +C\int_0^T\int_0^1\left(|\nb_x|^2v_x^2+u_x^2|\nb|^2+|\nw_x|^2\right)dxdt \\&\leq
C+\frac{1}{2}\int_0^T\int_0^1\frac{|\nb_{xx}|^2}{v^2}dxdt+C\int_0^T\max_{x\in[0,1]}|\nb_x|^2dt +\max_{(x,t)\in[0,1]\times[0,T]}|\nb|^2 \\&\leq C+\frac{3}{4}\int_0^T\int_0^1\frac{|\nb_{xx}|^2}{v^2}dxdt+C \int_0^T\int_0^1|\nb_x|^2dxdt \\&\quad  +C\sup_{0\le t\le T}\int_0^1|\nb|^2dx +\frac{1}{4}\sup_{0\le t\le T}\int_0^1|\nb_x|^2dx \\&\leq C +\frac{3}{4}\int_0^T\int_0^1\frac{|\nb_{xx}|^2}{v^2}dxdt+\frac{1}{4}\sup_{0\le t\le T}\int_0^1|\nb_x|^2dx,\ea\enn which implies
\bn\ba\la{zzz}\sup_{0\le t\le T}\int_0^1|\nb_x|^2dx+\int_0^T\int_0^1|\nb_{xx}|^2dxdt\leq C.\ea\en
Hence, \be\la{xxx} \max_{(x,t)\in [0,1]\times [0,T]}|\nb|^2\le C+C \sup_{0\le t\le T}\int_0^1|\nb_x|^2dx\le C.\ee

Finally, it follows from   \eqref{kkk}, \eqref{zzz}, \eqref{uiz1}, and  \eqref{xxx}  gives
\bnn\ba\la{rrr}\int_0^T\int_0^1|\nb_t|^2dxdt&\leq C\int_0^T\int_0^1\left(|\nb_{xx}|^2+|\nb_x|^2v_x^2+|\nw_x|^2+|\nb|^2u_x^2\right)dxdt\\&\leq C+C\int_0^T\left(\max_{x\in[0,1]}|\nb_x|^2+ \int_0^1u_x^2dx\right)dt\\&\leq C+C\int_0^T\int_0^1\left(|\nb_x|^2+|\nb_{xx}|^2\right)dxdt\\&\leq C.\ea\enn
 Combining this,  \eqref{ddd}, \eqref{ooo},  and  \eqref{zzz}  gives \eqref{aas1} and finishes the proof of Lemma \ref{lemmy0}. \thatsall

For further uses, we need the following estimate on the $L^2((0,1)\times(0,T))$-norm of $\te_x$ for $\b\in [0,1].$

 \begin{lemma}\la{lemm6} If $ \b\in[0,1],$ there exists a   positive constant $C $  such that for all $T> 0,$
  \be\ba\la{sq2} \int_0^T\int_0^1  \te_x^2dxdt
\le C+C \int_0^T\left(\int_0^1 u_x^2dx\right)^2dt .\ea\ee \end{lemma}
\pf    First, multiplying \eqref{2.4}   by $\te^{1-\frac{\b}{2}}$ and integration by parts gives
\be \la{p2}\ba  &\frac{2}{4-\b}\left(\int_0^1\te^{2-\frac{\b}{2}}dx\right)_t +\frac{(2-\b)}{2}\int_0^1 \frac{\te^{\frac{\b}{2}}\te_x^2}{v}dx\\& =-\int_0^1 \frac{\te^{2-\frac{\b}{2}} }{v}u_xdx+\int_0^1\frac{\te^{1-\frac{\b}{2}}( u_x^2+|\nb_x|^2+|\nw_x|^2)}{v} dx \\&=\int_0^1 \frac{\left(\bar\te^{2-\frac{\b}{2}}-\te^{2-\frac{\b}{2}}\right) }{v}u_xdx + \left(1-\bar\te^{2-\frac{\b}{2}} \right)\int_0^1 \frac{ u_x}{v}dx\\&\quad-\int_0^1 \frac{ u_x}{v}dx+\int_0^1\frac{\te^{1-\frac{\b}{2}}( u_x^2+|\nb_x|^2+|\nw_x|^2)}{v} dx\\ &\le C\int_0^1  \left|\te^{2-\frac{\b}{2}}-\bar\te^{2-\frac{\b}{2} } \right|\left| u_x\right|dx  +CV(t)-\left(\int_0^1 \ln v dx\right)_t\\&\quad + \int_0^1\frac{\te^{1-\frac{\b}{2}}( u_x^2+|\nb_x|^2+|\nw_x|^2)}{v} dx,\ea\ee where in the last inequality we have used \eqref{ux}.
Direct calculations yield that for any $\de>0$
\be\la{p4}\ba &\int_0^1  \left|\te^{2-\frac{\b}{2}}-\bar\te^{2-\frac{\b}{2} } \right|\left| u_x\right|dx \\&\le C\max_{x\in[0,1]}\left|\te^{1-\frac{\b}{4}}-\bar\te^{1-\frac{\b}{4}} \right|\left(\int_0^1\left(\te^{2-\frac{\b}{2}}+1\right)dx \right)^{1/2}\left(\int_0^1u_x^2dx\right)^{1/2}\\&\le \de\left(\int_0^1\te^{-\frac{\b}{4}}|\te_x|dx\right)^2+C (\de) \int_0^1\left(\te^{2-\frac{\b}{2}}+1\right)dx \int_0^1u_x^2dx  \\&\le C\de    \int_0^1\te^{\b/2} \te_x^2dx +C(\de)V(t)+C (\de) \int_0^1\left(\te^{2-\frac{\b}{2}}+1\right)dx \int_0^1u_x^2dx ,\ea\ee
 and that for any $\delta>0$\be\la{p5}\ba&\int_0^1\frac{\te^{1-\frac{\b}{2}} u_x^2}{v} dx \\ &\le C\left(\max_{x\in[0,1]}\left|\te^{ 1-\frac{\b}{2} }-\bar \te^{1-\frac{\b}{2} }\right|+1\right)   \int_0^1    u_x^2  dx  \\&\le C\int_0^1 \te^{-\frac{\b}{2}}|\te_x|dx\int_0^1   u_x^2  dx+C\int_0^1   u_x^2  dx  \\&\le \delta\int_0^1\left(\te^{\b-2} +\te^{\frac{\b}{2}}\right) \te_x^2dx+C(\delta)\left(\int_0^1   u_x^2  dx\right)^2+C\int_0^1   u_x^2 dx. \ea\ee
Putting \eqref{p4}  and \eqref{p5}  into \eqref{p2},   choosing $\de$ suitably small, and using \eqref{2.2},  \eqref{uiz1}, and the Gronwall inequality, one obtains \bnn\la{p1}\ba   \int_0^1\te^{2-\b/2}dx  +\int_0^T\int_0^1 \te^{\b/2}\te_x^2 dxdt\le C+C \int_0^T\left(\int_0^1 u_x^2dx\right)^2dt,\ea\enn which together with \eqref{2.2} implies
  \bnn\ba  \int_0^T\int_0^1  \te_x^2dxdt& \le C \int_0^T\int_0^1  \left(\te^{\b-2}+\te^{\b/2}\right) \te_x^2dxdt\\&\le C+C\int_0^T\int_0^1  \te^{\b/2} \te_x^2dxdt\\&
\le C+C \int_0^T\left(\int_0^1 u_x^2dx\right)^2dt .\ea\enn
 This gives \eqref{sq2} and finishes the proof of Lemma \ref{lemm6}.\thatsall

\begin{lemma}\la{lemmy11} For $\b\ge 0,$ there is a positive constant $C$ such that  for all $T> 0,$
\be \la{aaas1}\ba
&\sup_{0\le t\le T} \int_0^1 u_x^2 dx   +\int_0^T\int_0^1\left( u_t^2 +u_{xx}^2 +\te_x^2\right)dx dt\leq C.\ea\ee
\end{lemma}

\pf
First, rewriting \eqref{1.2} as
\be\ba\la{eq20} u_t=\frac{  u_{xx}}{v}- \frac{v_x   }{v^2 } u_x-\frac{\te_x}{v}+\frac{\te v_x}{v^2}-\nb\cdot\nb_x ,\ea\ee multiplying \eqref{eq20} by $u_{xx}$ and integrating the result over $(0,1) $,  we have
 \be\ba\la{eee} &\frac{1}{2}\frac{d}{dt} \int_0^1u_x^2dx+ \int_0^1\frac{  u_{xx}^2}{v}dx   \\&\leq  \frac{1}{2} \int_0^1\frac{  u_{xx}^2}{v}dx +C \int_0^1\left(\te_x^2+\te^2v_x^2+|\nb|^2|\nb_x|^2+u_x^2v_x^2 \right)dx   . \ea\ee
Direct computation yields that for any $\de>0,$
\be\la{qu1}\ba &\int_0^1\left(\te_x^2+\te^2v_x^2+|\nb|^2|\nb_x|^2+u_x^2v_x^2 \right)dx  \\ &\le C\left(\max_{x\in[0,1]} u_x^2+\max_{x\in[0,1]}\left(\te-\bar\te\right)^2+1\right)\int_0^1  v_x^2dx +C\int_0^1 \left(|\nb_x|^2 +  \te_x^2\right)dx\\ &\le C\max_{x\in[0,1]} u_x^2+C\max_{x\in[0,1]}\left(\te-\bar\te\right)^2  +C\int_0^1 ( v_x^2+|\nb_x|^2)dx +C\int_0^1 \te_x^2dx\\ &\le \delta \int_0^1 u_{xx}^2dx+ C(\delta) \int_0^1 u_x^2dx   +C\int_0^1  ( v_x^2+|\nb_x|^2)dx +C\int_0^1 \te_x^2dx,\ea\ee
where in the last inequality  we have used \eqref{ux2}.
Putting \eqref{qu1} into \eqref{eee} and choosing $\de$ suitably small yields
\be\ba\la{le5eq1}  \int_0^1 u_x^2dx +\int_0^T\int_0^1 u_{xx}^2 dxdt\le C   +C\int_0^T \int_0^1\te_x^2dx dt,\ea\ee due to    \eqref{aas1}  and \eqref{uiz1}.

Next, on the one hand, if $\b>1,$ choosing $p=\b-1$ in \eqref{ljj1} shows
\be\ba\la{sq1} \int_0^T\int_0^1  \te_x^2dxdt \le C, \ea\ee
 which  along with    \eqref{le5eq1}  gives
\be\ba\la{lm5w}  \sup_{0\le t\le T}\int_0^1 u_x^2dx+\int_0^T\int_0^1 u_{xx}^2dxdt+\int_0^T\int_0^1\te_x^2dxdt\le C.\ea\ee
 On the other hand, if $\b\in [0,1],$ it follows from \eqref{le5eq1}, \eqref{sq2},  \eqref{uiz1},      and Gronwall's inequality that \eqref{lm5w}  still holds.

Finally, it follows from \eqref{eq20}, \eqref{lm5w}, \eqref{qu1},  and \eqref{uiz1}  that
\bnn  \int_0^T\int_0^1 u_t^2dxdt\le C,\enn  which together with \eqref{lm5w} gives \eqref{aaas1} and finishes the proof of Lemma \ref{lemmy11}.\thatsall

\begin{lemma}\la{lemma70}For $\b\ge 0,$ there exists a positive constant $C$ such that for all $T>0,$\bn\ba\la{eq1} \sup_{0\le t\le T}\xix  \te_x^2dx+\int_0^T\xix \left( \te_t^2+\te_{xx}^2\right)dxdt\le C .  \ea\en
\end{lemma}

\pf
  First,
 multiplying \eqref{2.4} by $\te$ and integrating the result over $(0,1) $ yields
\be\ba\la{iii}&\frac{1}{2}\frac{d}{dt}\int_0^1\te^2dx+ \int_0^1\frac{\te^\b\te_x^2}{v}dx \\&= - \int_0^1\frac{\te^2-1}{v} u_x dx - \int_0^1\frac{u_x}{v}  dx+  \int_0^1\frac{\left( u_x^2 + |\nw_x|^2 + |\nb_x| ^2\right)\te}{v} dx  \\&\leq  C  \max_{x\in [0,1]}\left(u_x^2+ |\nw_x|^2 + |\nb_x| ^2+(\te-1)^2\right)  -\frac{d}{dt}\int_0^1\ln vdx , \ea\ee
 where we have used \eqref{wq1}. It follows from \eqref{vg4}, \eqref{aaas1}, and \eqref{uiz1} that
\bnn\ba \int_0^T \max_{x\in [0,1]}(\te-1)^2 dt& \le C \int_0^T \max_{x\in [0,1]}(\te-\bar\te)^2 dt+C \int_0^T  (1-\bar\te)^2 dt\\ & \le C \int_0^T \int_0^1\te_x^2dx dt+C \int_0^T  V(t) dt\le C,\ea\enn which together with \eqref{iii}, \eqref{aas1}, and  \eqref{aaas1} gives
\be \la{vc2}\int_0^T \int_0^1\te^{\b}\te_x^2dxdt\le C.\ee

Next, noticing that integration by parts leads to
\bnn\ba\la{eq3}  \xix \te^{\b}\te_t\left(\frac{\te^\b\theta_{x}}{v}\right)_{x}dx & =-\xix \frac{\te^\b\theta_{x}}{v}\left(\te^{\b}\te_t\right)_{x}dx\\& =-\xix \frac{\te^\b\theta_{x}}{v}\left(\te^{\b}\te_x\right)_{t}dx \\&=-\frac{1}{2}
\xix \frac{\left((\te^\b\theta_{x})^2\right)_t}{v}dx \\&=-\frac{1}{2} \left(\xix \frac{(\te^\b\theta_{x})^2}{v}dx\right)_t
-\frac{1}{2}\xix \frac{(\te^\b\theta_{x})^2u_x}{v^2}dx  ,\ea\enn
multiplying \eqref{2.4} by
$ \te^{\b}\te_t$ and integrating the resultant equality over (0,1), we have
\be\ba\la{eq4} & \xix  \te^\b\te_t^2dx+ \frac{1}{2} \left(\xix \frac{(\te^\b\theta_{x})^2}{v}dx\right)_t\\&=- \frac{1}{2}\xix \frac{(\te^\b\theta_{x})^2u_x}{v^2}dx+\xix  \frac{ \te^{\b }\theta_t\left(-\te u_x+  u_x^2+|\nw_x|^2+|\nb_x|^2\right)}{v }dx\\&\le C\max_{x\in[0,1]} (|u_x|\te^{\b/2})\xix \te^{3\b/2}\te_x^2dx +\frac{1}{2}\xix  \te^\b\te_t^2dx+C\xix  \te^{\b+2}u_x^2dx\\&\quad+C\xix \te^\b\left(u_x^4+|\nw_x|^4+|\nb_x|^4\right)dx%\\&\le C \xix  \te^{2\b}\te_x^2dx\xix  \te^\b\te_x^2dx  +\frac{1}{2}\xix  \te^\b\te_t^2dx\\&\quad+C \max_{x\in[0,1]} (\te^{\b+2}+ \te^\b \left(u_x^2+|\nw_x|^2+|\nb_x|^2\right))+C
\\&\le  C \xix  \te^{2\b}\te_x^2dx\xix  \te^\b\te_x^2dx+\frac{1}{2}\xix  \te^\b\te_t^2dx+\max_{x\in[0,1]}(\te^{1+3\b/4}-\bar \te^{1+3\b/4})^4\\&\quad+C \left(1+\int_0^1 \te^{2\b}\te_x^2dx \right)\max_{x\in[0,1]}\left(   u_x^2+  u_x^4 +|\nw_x|^4+|\nb_x|^4 \right) ,\ea\ee due to  \be \la{vc1} \max_{x\in[0,1]} (\te^{\b+1}-\bar\te^{ \b+1})^2\le C\int_0^1 \te^{2\b}\te_x^2dx.\ee

Next, combining \eqref{ux2} and \eqref{aaas1} gives
 \be\ba\la{eq5} \int_0^T\max_{x\in[0,1]}u_x^4dt&\leq C\int_0^T\int_0^1u_x^2dx\int_0^1u_{xx}^2dxdt \\ &\leq C\int_0^T \int_0^1u_{xx}^2dxdt\le C  .\ea\ee
  Using \eqref{aas1} and applying similar arguments to $\nb$ and $\nw$ implies
 \bn\la{eq6}\ba\int_0^T\max_{x\in[0,1]}|\nb_x|^4dt\leq C,\quad \int_0^T\max_{x\in[0,1]}|\nw_x|^4dt\leq C. \ea\en
Noticing that \be\la{eq9}\ba\max_{x\in[0,1]}\left(   \te^{1+3\b/4}-\bar \te^{1+3\b/4}\right)^4 &\le C\left(\xix  \te^{3\b/4}|\te_x|dx\right)^4\\ &\le   C \xix  \te^{2\b}\te_x^2dx\xix  \te^\b\te_x^2dx ,\ea\ee we   deduce from \eqref{eq4}, \eqref{vc2},   \eqref{eq5}--\eqref{eq9}, and the Gronwall inequality that
\be\ba\la{eq8} \sup_{0 \le t\le T}\xix  \left(\te^\b\theta_{x}\right)^2 dx+\int_0^T\xix  \te^\b\te_t^2dxdt\le C, \ea\ee
which together with \eqref{vc1} in particular gives
\be\la{eq10}\max_{(x,t)\in[0,1]\times[0,T]}\te(x,t)\le C.\ee

Next, it follows from \eqref{aaas1} and \eqref{eq10} that
\be\ba \la{lm8eq5} \int_0^T \int_0^1 \left(\te^{\b+1}- \bar \te^{\b+1}\right)^2dxdt &\le C\int_0^T \left(\int_0^1  \te^{{\b}}|\te_x| dx\right)^2dt
\\ &\le C\int_0^T \int_0^1   \te_x^2 dx dt   \le C,\ea\ee
which together with  \eqref{eq10},   \eqref{eq8}, and \eqref{2.4} shows
\be\ba \la{lm8rq5} &\int_0^T \left|\frac{d}{dt}\int_0^1 \left(\te^{\b+1}- \bar \te^{\b+1}\right)^2dx\right|dt \\&\le C\int_0^T \int_0^1 \left(\te^{\b+1}- \bar \te^{\b+1}\right)^2dxdt+C\int_0^T \int_0^1\left(\te^\b\te_t^2+ \bar\te_t^2\right)dxdt\\&\le C+C\int_0^T\int_0^1\left(u_x^2+|\nb_x|^2+ |\nw_x|^2\right)dxdt\\&\le C.\ea\ee

Thus, both \eqref{lm8eq5} and \eqref{lm8rq5} lead   to
\bnn \lim_{t\rightarrow \infty}\int_0^1 \left(\te^{\b+1}- \bar \te^{\b+1}\right)^2dx=0,\enn which  together with  \eqref{eq8} gives
\be \la{vc4}\ba  \max_{x\in[0,1]} (\te^{\b+1}-\bar\te^{ \b+1})^4&\le C\int_0^1 \left(\te^{\b+1}- \bar \te^{\b+1}\right)^2dx\int_0^1 \te^{2\b}\te_x^2dx\\&\le C\int_0^1 \left(\te^{\b+1}- \bar \te^{\b+1}\right)^2dx\rightarrow 0 \mbox{ as }t\rightarrow \infty.\ea\ee

It thus follows from \eqref{vc4}  and \eqref{qwe1}  that there exists some  $  T_0>0$ such that
\be\la{x5}\ba \te (x,t)\ge \al_1/2, \ea\ee   for all $(x,t)\in [0,1]\times [ T_0,\infty).$ Moreover, it follows from  \cite{ka1,hss1} that  there exists some constant $C \ge  2/\al_1$
 such that  \bnn \theta(x,t)\ge  C^{-1} ,\enn  for all $(x,t)\in [0,1]\times [0,T_0].$ Combining this, \eqref{x5}, and \eqref{eq10} yields that  for all $(x,t)\in[0,1]\times[0,\infty),$  \be\ba\la{tq2}  C^{-1}\le \te(x,t)\le C  , \ea\ee
which together with \eqref{eq8} gives
\be\ba\la{eq11}\sup_{0 \le t\le T}\xix  \theta_{x}^2 dx+\int_0^T\xix  \te_t^2dxdt\le C. \ea\ee

Finally, it follows from \eqref{2.4} that
\bnn\ba \frac{\te^\b\te_{xx}}{v}= -\frac{\b \te^{\b-1}\te_x^2}{v}+\frac{\te^\b\te_x v_x}{v^2}- \frac{  u_x^2+|\nb_x|^2+|\nw_x|^2}{v}+ \frac{ \te u_x}{v}+\te_t,\ea\enn
which together with \eqref{tq2}, \eqref{uiz1},  \eqref{eq5}, \eqref{eq6}, \eqref{aaas1},   \eqref{eq11}, and \eqref{ux2} yields
\bnn\ba\la{eq12}\int_0^T\xix \te_{xx}^2dxdt  & \le C\int_0^T\xix \left(\te_x^4+\te_x^2v_x^2+u_x^4+|\nb_x|^4+|\nw_x|^4+u_x^2+\te_t^2\right)dxdt\\ & \le C+ C\int_0^T\max_{x\in[0,1]}\te_x^2 dt\\ & \le C+   \frac12\int_0^T  \xix \te_{xx}^2 dx   dt.\ea\enn
 Combining this with \eqref{eq11}  proves \eqref{eq1}
and  finishes the proof of \lemref{lemma70}.\thatsall

Finally, we have the following nonlinearly exponential stability of the strong solutions.
\begin{lemma} There exist  some positive constants $C$ and $\eta_0 $ both
 depending only on $\b,\|(v_0,u_0,\theta_0,\nb_0,\nw_0 )\|_{H^1(0,1)},
 \inf\limits_{x\in [0,1]}v_0(x),$ and $ \inf\limits_{x\in [0,1]}\theta_0(x) $ such that
\be \la{klq1} \|(v-1,u,\te-1,\nb,\nw)(\cdot,t)\|_{H^1(0,1)}\le Ce^{-\eta_0t}.\ee\end{lemma}

\pf Noticing that all the constants $C$ in Lemmas \ref{lemq1}, \ref{lemmy0}, \ref{lemmy11}, and \ref{lemma70}    are independent of $T,$ we have
\be\la{zma1} \int_0^\infty \left| \frac{d}{dt}\left(\|v_x \|_{L^2}^2, \|u_x \|_{L^2}^2, \|\te_x \|_{L^2}^2, \|\nb_x \|_{L^2}^2 ,\|\nw_x \|_{L^2}^2\right) \right|dt\le C,\ee where we have used \bnn
 \int_0^1u_xu_{xt}dx=- \int_0^1u_tu_{xx}dx.\enn It thus follows from \eqref{zma1}, \eqref{uiz1}, and  \eqref{eq1} that
\bnn \lim_{t\rightarrow \infty}\|(v_x,u_x,\te_x,\nb_x,\nw_x)(\cdot,t)\|_{L^2(0,1)}=0,\enn which in particular implies
\be \la{bd1} \lim_{t\rightarrow \infty}\|(v-1,u ,\te-1,\nb,\nw)(\cdot,t)\|_{H^1(0,1)}=0.\ee
With \eqref{bd1} at hand, the proof of \eqref{klq1} is standard   (c.f. \cite{ok}). \thatsall

%\newpage

 \end{document}